\documentclass[a4paper,12pt]{amsart}
\usepackage{color}
\usepackage{amsmath,amssymb,amscd,amsfonts,amsthm,stmaryrd}
\usepackage[mathscr]{eucal}
\usepackage{xcolor}
\usepackage{graphicx}
\usepackage{tikz-cd}
\usepackage{adjustbox}

\usepackage{nicefrac}

\usepackage{color}

\numberwithin{equation}{section}

\def\rank{\operatorname{rank}}

\def\rad{\operatorname{rad}}

\def\acts{\curvearrowright}
\def\D{\partial}

\def\R{\mathbb R}
\def\Z{\mathbb Z}

\def\N{\mathbb N}

\def\al{\alpha}

\def\ka{\kappa}

\def\eps{\epsilon}
\def\ga{\gamma}
\def\Ga{\Gamma}

\def\la{\lambda}

\def\om{\omega}

\def\si{\sigma}
\def\Si{\Sigma}

\def\tits{\partial_{T}}
\def\geo{\partial_{\infty}}

\def\wlim{\mathop{\hbox{$\om$-lim}}}

\def\acts{\curvearrowright}

\def\Isom{\mathop{\hbox{Isom}}}

\def\Ant{\operatorname{Ant}}
\def\Ce{\operatorname{Centers}}

\def\<{\langle}
\def\>{\rangle}

\theoremstyle{plain}
\newtheorem{thm}{Theorem}[section]
\newtheorem{lem}[thm]{Lemma}
\newtheorem{prop}[thm]{Proposition}
\newtheorem{cor}[thm]{Corollary}
\newtheorem{slem}[thm]{Sublemma}

\newtheorem{introthm}{Theorem}

\newtheoremstyle{named}{}{}{\itshape}{}{\bfseries}{.}{.5em}{\thmnote{#3} #1}
\theoremstyle{named}

\theoremstyle{definition}
\newtheorem{dfn}[thm]{Definition}

\theoremstyle{remark}

\newcommand{\bcl}{\begin{claim}}
\newcommand{\ecl}{\end{claim}}
\newcommand{\bcor}{\begin{cor}}
\newcommand{\ecor}{\end{cor}}
\newcommand{\bdfn}{\begin{dfn}}
\newcommand{\edfn}{\end{dfn}}
\newcommand{\ben}{\begin{enumerate}}
\newcommand{\bit}{\begin{itemize}}
\newcommand{\blem}{\begin{lem}}
\newcommand{\bslem}{\begin{slem}}
\newcommand{\bprop}{\begin{prop}}
\newcommand{\bthm}{\begin{thm}}
\newcommand{\een}{\end{enumerate}}
\newcommand{\eit}{\end{itemize}}
\newcommand{\elem}{\end{lem}}
\newcommand{\eslem}{\end{slem}}
\newcommand{\eprop}{\end{prop}}
\newcommand{\ethm}{\end{thm}}

\usepackage{hyperref}

\begin{document}

\title[Rigidity of ideal symmetric sets]{Rigidity of ideal symmetric sets}

\author{Stephan Stadler}

\newcommand{\Addresses}{{\bigskip\footnotesize
\noindent Stephan Stadler,
\par\nopagebreak\noindent\textsc{Max Planck Institute for Mathematics, Vivatsgasse 7, 53111 Bonn, Germany}
\par\nopagebreak
\noindent\textit{Email}: \texttt{stadler@mpim-bonn.mpg.de}

}}




\begin{abstract}
We characterize higher rank model geometries -- Riemannian symmetric spaces, Euclidean buildings
and products -- among Hadamard spaces by  using antipodal sets at infinity.
\end{abstract}

\maketitle

\section{Introduction}
\subsection{The main result}
A {\em Hadamard space} $X$ is a locally compact and geodesically complete CAT(0) space.
A subset  of its ideal boundary $A\subset\geo X$ is called {\em symmetric},
if every geodesic line in $X$ has either both or none of its ideal endpoints in $A$.
If $M$ is a smooth Hadamard manifold whose ideal boundary $\geo M$ contains a non-trivial closed symmetric subset,
then Eberlein  has shown that the holonomy group of $M$ cannot be transitive \cite{E_symII}. It follows from Berger's holonomy theorem 
 that $M$ is either a higher rank symmetric space or else splits as a non-trivial product. 

\begin{introthm}\label{thm_mainA}
Let $X$ be a cocompact Hadamard space. 
Suppose that $\geo X$ contains a non-trivial closed symmetric subset.
Then $X$ is a Riemannian symmetric space, a Euclidean building or splits as a non-trivial product.
\end{introthm}

Ideal symmetric subsets as above can often be constructed if the space has enough symmetries. We will show
that if the isometry group of a Hadamard  space $X$ is  large in a dynamical sense, then either there  
exists an ideal closed symmetric subset or $X$ resembles a space of negative curvature.
  
\subsection{Applications to rigidity}
	
Let $X$ be a Hadamard space and $\Ga$ a group of isometries of $X$.
Then $\Ga$ is said to satisfy the {\em duality condition}, if the set of $\Ga$-recurrent geodesics is dense. (See Section~\ref{sec_dyn}.) 
Recall that if a group $\Ga$ acts on a Hadamard manifold isometrically, properly discontinuously and with finite covolume, then $\Ga$
satisfies the duality condition \cite[p.~39]{ballmannbook}.

\begin{introthm}\label{thm_mainB}
Let $X$ be a cocompact Hadamard space. Suppose that $\Ga$ is a group of isometries satisfying the duality condition.
Then either the geodesic flow of $X$ contains a dense orbit modulo $\Ga$, or else  
$X$ is a Riemannian symmetric space, a Euclidean building or splits as a non-trivial product.
\end{introthm}

We can strengthen the conclusion for Hadamard spaces whose Tits boundaries are not everywhere locally non-compact.
Recall that a complete geodesic in a Hadamard space has {\em rank 1}, if it does not bound a flat half-plane \cite{ballmann_diss}.

\begin{introthm}\label{thm_mainC}
Let $X$ be a cocompact Hadamard space. Let $\Ga$ be a group of isometries of $X$ satisfying the duality condition.
Suppose that the Tits boundary contains a relatively compact open subset.
Then either $X$ contains a $\Ga$-periodic rank 1 geodesic, or else  
$X$ is a Riemannian symmetric space, a Euclidean building or splits as a non-trivial product.
\end{introthm}

If $X$ contains a $\Ga$-periodic rank 1 geodesic, then $X$ and $\Ga$ display many hyperbolic features \cite[Theorem~B]{ballmannbook}. 
By structure theory of locally compact spaces with upper curvature bounds \cite{LN_gcba}, the Tits boundary contains a relatively compact open set if
and only if it contains a regular point. Recall that a periodic $n$-flat in $X$ is called {\em Morse}, if it does not bound a flat $(n+1)$-half-space. 
It follows from our earlier work \cite{St_rrI} that a periodic Morse flat leads to regular points in the Tits boundary.
The existence of a periodic Morse flat can be attributed to the acting group. For an arbitrary finitely generated discrete group $\Ga$ we introduce its {\em Morse rank} 
\[\rank_{Mo}(\Ga)\in\N\cup\{\infty\}.\]
If such a group acts geometrically on a CAT(0) space $X$, then $\rank_{Mo}(\Ga)$ is finite if and only if $X$ contains a $\Ga$-periodic Morse flat. (See Section~\ref{sec_morse}.) 
If a group $\Ga$ acts geometrically on a Riemannian symmetric space of noncompact type or more generally on a Hadamard manifold $M$, then 
$\rank_{Mo}(\Ga)=\rank(M)$ holds.

Recall that a Hadamard space $X$ is called {\em rank 1}, if it contains a rank 1 geodesic which is an axis for some isometry of $X$. 

\begin{introthm}\label{thm_mainD}
Let $X$ be a Hadamard space with a geometric group action $\Ga\acts X$. Suppose that $\Ga$ satisfies the duality condition and
has finite Morse rank.
Then $X$ has a unique decomposition into irreducible factors
\[X\cong\R^j\times X_1\times\ldots\times X_k\times Y_1\times\ldots\times Y_l\times M_1\times\ldots\times M_m\]
where $X_a$ are rank 1 spaces, $Y_b$ are irreducible Euclidean buildings of higher rank and $M_c$
are irreducible Riemannian symmetric spaces of higher rank.
\end{introthm}

This has consequences for the algebraic structure of the acting groups:

\begin{introthm}[{Tits Alternative}]\label{thm_mainE}
Let $X$ be a Hadamard space with a geometric group action $\Ga\acts X$. Suppose that $\Ga$ satisfies the duality condition and
has finite Morse rank. Then either $\Ga$ is virtually abelian and $X$ is flat or $\Ga$ contains a non-abelian free subgroup.
\end{introthm}

\subsection{On the proof of Theorem~\ref{thm_mainA}}

Besides its relevance for Rank Rigidity, our main theorem is motivated by Lytchak's rigidity result for spherical buildings and joins.
To state the theorem, recall that a subset $A$ in a CAT(1) space is called {\em symmetric}, if every point at distance at least $\pi$ from a point in $A$ 
has to lie in $A$. Further, a CAT(1) space is called {\em geodesically complete}, if every geodesic segment extends to a complete local geodesic.

\bthm[{\cite[Main Theorem]{Ly_rigidity}}]
Let $Z$ be a finite-dimensional geodesically complete CAT(1) space. If $Z$ contains a non-trivial closed symmetric subset,
then $Z$ is a spherical building or a join. 
\ethm

The Tits boundary of a cocompact Hadamard space is a finite-dimensional CAT(1) space. 
However, it is rarely ever geodesically complete, even if its diameter is $\pi$. (See the discussion in \cite[p.~711]{St_rrII}.)
As in \cite{St_rrII}, the proof of Theorem~\ref{thm_mainA} follows the strategy of \cite[Main Theorem]{Ly_rigidity} and studies the submetry $\delta:\tits X\to\Delta$
induced by minimal closed symmetric sets. The major novelty of this paper over \cite{St_rrII} is that we do not need any additional assumptions (like the presence of a periodic Morse flat)
to understand the structure of the submetry $\delta$. Recall that a point in a metric space is called {\em regular}, if it has a neighborhood homeomorphic to an
open set in a Euclidean space. We use a completely new idea to show the following dichotomy for a point $\xi$ in a top-dimensional round sphere 
$\hat\si\subset\tits X$. 
Either $\xi$ is a regular point or there exists a subsphere $\si\subset\hat\si$ through the point $\xi$ such that the restricted submetry 
$\delta|_{\hat\si}$ is invariant under reflection along $\si$. So either $\hat\si$ contains many regular points or else $\delta|_{\hat\si}$
is very symmetric -- a win-win situation.

\subsection{Organization}
In Section~\ref{sec_pre} we recall basics and agree on notation. In Section~\ref{sec_sub}
we discuss ideal symmetric sets and first properties of their induced submetry $\delta$ on the Tits boundary.
We show that $\delta$ restricts to a submetry with the same base on every top-dimensional round sphere.
 We also introduce the notion of singular spheres
and explain how they lead to symmetries of restricted submetries. In Section~\ref{sec_ind_sub} we produce many singular spheres
using the isometry group. In Section~\ref{sec_irr} we deduce splitting results from certain configurations of non-branching piecewise horizontal geodesics.
We show that for an irreducible space with a non-trivial submetry $\delta$ the Tits boundary has to be geodesically complete. In Section~\ref{sec_main}
we provide proofs of the main result and its applications. The final Section~\ref{sec_tits} is devoted to the Tits Alternative.

\subsection{Acknowledgments}
I want to thank Alexander Lytchak for reading a preliminary version of this paper and for his useful comments.
I was supported by DFG grant SPP 2026.

\section{Preliminaries}\label{sec_pre}

References for this section are \cite{AKP, ballmann_geo, BH, KleinerLeeb, Ly_rigidity}.

\subsection{Metric spaces}
Euclidean $n$-space with its flat metric will be denoted by $\R^n$ and the unit sphere by $S^{n-1}\subset\R^n$.
We denote  the distance between two points $x$ and $y$ in a metric space $X$ by $|x,y|$. 
For $x\in X$ and $r>0$, we denote by $B_r(x)$ and $\bar B_r(x)$ the open and closed $r$-ball around $x$, respectively.
If $A\subset X$ denotes a subset, then $|x,A|$ and  refers to the greatest lower bound for distances from points in $A$ to $x$.
Similarly, $|x,A|_H$ refers to the {\em Hausdorff distance} -- smallest upper bound for distances from points in $A$ to $x$.
More generally, for subsets $A, A'\subset X$ we denote by $|A,A'|$ and $|A,A'|_H$ the greatest lower bound for distances from points in $A$ to points in $A'$ and
the Hausdroff distance, respectively.
The {\em radius}  of a set $A \subset X$ is defined by 
\[\rad_X(A):=\inf_{x\in X}|x,A|_H.\]
The set of {\em centers} of $A$ in $X$ is given by
\[\Ce_X(A):=\{x\in X|\ |x,A|_H=\rad_X(A)\}.\]
However, if the respective space is clear from the context, then we will suppress the specification.
Notice that $\rad_X (A)$ and $\Ce_X (A)$ depend on the ambient space $X$.
$\Ce(A)$ is a closed subset which is preserved by every isometry which preserves the set $A$.

A \emph{geodesic}
is an isometric embedding of an interval. It is called a {\em geodesic segment}, if it is compact.
The {\em endpoints} of a geodesic segment $c$ are denoted by $\D c$.
A geodesic segment $c$ {\em branches} at an endpoint $y\in \D c=\{x,y\}$, if there are geodesics $c^\pm$ starting in $x$ which  strictly contain $c$
and such that $c^-\cap c^+=c$.

 A \emph{triangle} is a union of three geodesics connecting three points.
If $x,y,z$ are point with unique geodesics between them, then we denote the corresponding triangle by $\triangle(x,y,z)$.
$X$  is a \emph{geodesic metric space} if
any pair of   points of $X$
is connected by a geodesic.
It is \emph{geodesically complete} if every geodesic segment is contained in a complete local geodesic.

\subsection{Submetries}

In this section we recall the notion of a submetry. Originally the concept was introduced by Berestovskii in \cite{Ber_sub}
as a metric analog of a Riemannian submersion. If $B$ is a spherical building, then the natural $1$-Lipschitz map $\delta: B\to\Delta_{mod}$ 
which folds the building onto is model chamber is a submetry. Vice versa,  CAT(1) spaces
which admit certain kinds of submetries are spherical buildings \cite{Ly_rigidity}.
This explains the relevance of submetries for us. 
\medskip

Two subsets $X^\pm$ in a metric space $X$ are called {\em equidistant}, if for every point $x^\pm\in X^\pm$
there exists a point $y^\mp\in X^\mp$ such that $|x^\pm,y^\mp|=|X^\pm,X^\mp|$ holds.

A decomposition of a metric space $X$ into closed pairwise equidistant subsets $\{X_y\subset X|\ y\in Y\}$ induces a natural metric on the quotient space $Y$:
\[|y,y'|:=|X_y,X_{y'}|.\]
The natural projection $\delta:X\to Y$ is a submetry in the sense of \cite{Ber_sub}.
\bdfn
A {\em submetry} $\delta:X\to Y$ between metric spaces is a map which sends for every point $x$ in $X$ and each radius $r\geq 0$
the closed $r$-ball around $x$ surjectively onto the closed $r$-ball around $\delta(x)$.
\edfn
Submetries on a metric space $X$ are in one-to-one correspondence to closed equidistant decompositions of $X$.
Each submetry is a 1-Lipschitz map.

If $\delta:X\to Y$ is a submetry, then
two points $x,x'$ in $X$ are called $\delta$-near, if $|x,x'|=|\delta(x),\delta(x')|$ holds.
In this case, $\delta$ maps any geodesic between $x$ and $x'$ isometrically onto a geodesic between $\delta(x)$
and $\delta(x')$. Such geodesics in $X$ are called {\em horizontal}. More generally, a geodesic in $X$ is called {\em piecewise horizontal},
if it is the union of finitely many horizontal geodesics. 

We collect a couple of results on submetries of round spheres which will be used later on.
Most of these results hold in greater generality but we restrict to the cases of interest for us.

Let $\delta:S^n\to\Delta$ be a submetry. Then $\Delta$ is a compact Alexandrov space of dimension $k\leq n$ and curvature at least 1 \cite[Proposition~3.1]{KaLy_sub},
\cite{BGP}.

It allows for a natural disjoint decomposition into strata $\Delta^l\subset\Delta$, $0\leq l\leq k$.
The $l$-dimensional stratum $\Delta^l$ is
the set of points $y\in\Delta$ whose tangent space $T_y\Delta$ splits off  $\R^l$ as a direct
factor but not $\R^{l+1}$.

\bthm[{\cite[Theorem~1.6]{KaLy_sub}}]
Let $\delta:S^n\to\Delta$ be a submetry with $\dim(\Delta)=k$. For any $0\leq l\leq k$, the stratum $\Delta^l$ is
an $l$-dimensional topological manifold which is locally closed and locally
convex in $\Delta$. The top-dimensional stratum $\Delta^k$ is open and globally convex.
\ethm

The top-dimensional stratum is sometimes referred to as the set of {\em regular points},
$\Delta_{reg}:=\Delta^k$. It is an open dense convex subset in $\Delta$ of full Hausdorff $k$-measure \cite{BGP}. 

\bthm[{\cite[Theorem~10.5]{KaLy_sub}}]\label{thm_submersion}
Let $\delta:S^n\to\Delta$ be a submetry. Then $\delta:\delta^{-1}(\Delta_{reg})\to\Delta_{reg}$
is a $\mathcal{C}^{1,1}$ Riemannian submersion. 
\ethm

We call a submetry $\delta:S^n\to\Delta$ {\em transnormal}, if any extension of a $\delta$-horizontal geodesic in $S^n$ is a piecewise horizontal
geodesic.
This is equivalent to saying that an antipode of a direction tangent to a $\delta$-horizontal geodesic is again tangent to a $\delta$-horizontal geodesic.
Since the horizontal space $H_x\subset\Si_x S^n$ is convex for any submetry \cite[Proposition~12.5]{KaLy_sub}, it follows that for transnormal submetries
horizontal spaces are round spheres. We even have:

\bprop[{\cite[Proposition~12.5]{KaLy_sub}}]\label{prop_C11}
Let $\delta:S^n\to\Delta$ be a transnormal submetry. Then every fiber of $\delta$ is a $\mathcal{C}^{1,1}$ submanifold in $S^n$.
\eprop

The following is an easy consequence of \cite[Theorem~1.2]{La_orbi}.

\blem[{\cite[Lemma~2.5]{St_rrII}}]\label{lem_polytope}
Let $\delta:S^n\to\Delta$ be a submetry such that every geodesic in $S^n$ is piecewise $\delta$-horizontal. Then
$\Delta$ is a spherical orbifold and $\delta$ is a covering of Riemannian orbifolds.
\elem

Let $\hat\delta:\hat\si\to\Delta$ be a transnormal submetry of a round $n$-sphere $\hat\si$.
A round subsphere $\si\subset\hat\si$ is called a {\em reflecting sphere (for $\hat\delta$)}, if
$\hat\delta$ is invariant under reflection along $\si$. More precisely, if $\si\subset\hat\si$ is  a hypersphere,
then we ask that $\hat\delta$ is invariant under reflection along $\si$. If $\si\subset\hat\si$ has arbitrary codimension, we ask that $\delta$
is invariant under reflection along every hypersphere $\si'\subset\hat\si$ which contains $\si$.

The existence of reflecting spheres for a submetry $\hat\delta:\hat\si\to\Delta$ is rather restrictive.
For instance, if there is a single reflecting sphere, then  the  boundary of $\Delta$ is non-trivial.

\bprop\label{prop_splitsub}
Let $\hat\delta:\hat\si\to\Delta$ be a transnormal submetry of a round $n$-sphere $\hat\si$. Let  $\si_1,\ldots,\si_j\subset\hat\si$ with $\dim(\si_i)=n-d_i$ 
be reflecting spheres. Suppose that $\sum_{i=1}^j d_j=n+1$ and  $\bigcap_{i=1}^{j}\si_i=\emptyset$. 
Then either every geodesic in $\hat\si$  is piecewise $\hat\delta$-horizontal,
or there exists a disjoint decomposition $I^-\dot\cup I^+=\{1,\ldots,j\}$ such that $\hat\si=\si^-\circ\si^+$
where $\si^\pm=\bigcap_{i\in I^\pm}\si_i$ and every geodesic $\xi^-\xi^+$ with $\xi^\pm\in\si^\pm$ is piecewise $\hat\delta$-horizontal. 
\eprop

\proof
Choose for every $\si_i$ hyperspheres $\si_{ik}$, $1\leq k\leq d_i$,
such that $\si_i=\bigcap_{k=1}^{d_i}\si_{ik}$. By assumption, every hypersphere $\si_{ik}$ is a reflecting sphere for $\hat\delta$.
Thus the claim  follows from \cite[Proposition~2.9]{St_rrII}.
\qed

\subsection{Spaces with an upper curvature bound}

For $\kappa \in \R $, let $D_{\kappa}  \in (0,\infty] $ be the diameter of the  complete, simply connected  surface $M^2_{\kappa}$
of constant curvature $\kappa$. A complete  metric space is called a \emph{CAT($\kappa$) space}
if any pair of its points with distance less than $D_{\kappa}$ is connected by a geodesic and if
 all triangles
with perimeter less than $2D_{\kappa}$
are not thicker than
 the \emph{comparison triangle} in $M_{\kappa} ^2$. In particular, geodesics between points of distance less than $D_{\kappa}$
are unique. 

For any CAT($\kappa$) space  $X$, 
 and points $x\neq y$ at distance less than $D_{\kappa}$, we denote the unique geodesic
 between $x$ and $y$ by $xy$.   For $y,z \neq x$, the angle at $x$ between $xy$ and $xz$
will be denoted by $\angle_x(y,z)$.
The \emph{space of directions} or {\em link}
 at a point  $x\in X$ is denoted by $\Si_x X$, it is a CAT(1) space when equipped with the angle metric.
Its elements are called {\em directions}.
A direction $v$ is called {\em genuine},
if there exists a geodesic starting in direction $v$.

\subsection{Dimension}
 A natural notion of dimension $\dim (X)$ for a CAT($\ka$) space $X$ was
  introduced by
Kleiner in \cite{Kleiner}, originally referred to as {\em geometric dimension}.
It vanishes precisely when the space is discrete.
In general, it is defined inductively:
\[\dim (X)= \sup _{x\in X} \{ \dim (\Sigma _xX) +1 \}.\]
For instance, in a 1-dimensional CAT($\ka$) space every link is discrete. 
The geometric dimension coincides with the supremum of 
topological dimensions\footnote{Here topological dimension corresponds to Lebesgue covering dimension.} of compact subsets in $X$ \cite{Kleiner}. 
If the dimension of $X$ is finite, then $\dim(X)$ agrees with the largest number $k$
such that $X$ admits a bilipschitz embedding of an open subset in $\R^k$ \cite{Kleiner}.
The dimension of a locally compact and geodesically complete space is finite and agrees with the topological dimension 
as well as the Hausdorff dimension \cite{OT_cba, LN_gcba}.  

\subsection{CAT(1) spaces}

CAT(1) spaces play a particular role in our investigations as they appear as Tits boundaries and links of CAT(0) spaces.
Recall that the spherical join $Z_1\circ Z_2$ of two CAT(1) spaces $Z_1$ and $Z_2$ is a CAT(1) space of diameter $\pi$.

 A subset $C$ in a CAT(1) space $Z$ is called {\em $\pi$-convex}
if for any pair of points $x, y\in Z$ at distance less than $\pi$ the unique geodesic $xy$ is contained in $C$.
If $C$ is closed, then it is CAT(1) with respect to the induced metric.
For instance, a ball of radius at most $\frac{\pi}{2}$ is $\pi$-convex.

A subset  in a CAT(1) space is called {\em spherical}, if it 
embeds isometrically into a round sphere. A {\em (spherical) $n$-lune of angle $\theta\in[0,\pi]$} in a CAT(1) space is a closed convex subset $\la$ isometric to 
$S^{n-1}\circ[0,\theta]$. By the Lune Lemma \cite[Lemma~2.5]{BB_diam}, a geodesic bigon in a CAT(1)
space spans a $2$-lune. 
We will use the following more general result which follows immediately from \cite[Lemma~2.5]{BB_diam} and \cite[Lemma~4.1]{Ly_rigidity}.

\blem[Lune Lemma]\label{lem_lune}
Let $\xi$ and $\hat\xi$ be antipodes at distance $\pi$ in a CAT(1) space $Z$.
If $\ga^\pm$ are geodesics from $\xi$ to $\hat\xi$, then
$\ga^-\cup\ga^+$ is a spherical subset in $Z$. 
More precisely, if $\angle_\xi(\ga^-,\ga^+)\geq \pi$, then $\ga^-\cup\ga^+$ is a round $1$-sphere in $Z$;
and if  $\angle_\xi(\ga^-,\ga^+)=\theta< \pi$, then $\ga^-\cup\ga^+$ bounds a spherical lune of angle $\theta$.
Similarly, let $\tau^\pm\subset Z$ be round $n$-hemispheres with $\D\tau^-=\D\tau^+$ and let $\zeta^\pm\in Z$
denote their centers. If $|\zeta^-,\zeta^+|\geq\pi$, then $\tau^-\cup\tau^+$ is a round $n$-sphere in $Z$;
and if $|\zeta^-,\zeta^+|=\theta<\pi$, then $\tau^-\cup\tau^+$ bounds a spherical $(n+1)$-lune of angle $\theta$.
\elem

Two points in a CAT(1) space are called {\em antipodes}, if their distance is at least $\pi$.
A point in a CAT(1) space is called {\em regular}, if it has a neighborhood homeomorphic to an open subset in a Euclidean space.
All other points are called {\em singular}. 
A regular point is called
{\em $k$-spherical}, if it has a neighborhood isometric to an open set in the round $k$-sphere.

\subsection{CAT(0) spaces}

The {\em ideal boundary} of a CAT(0) space $X$, equipped with the cone topology, is denoted by $\geo X$. 
If $X$ is locally compact, then $\geo X$ is compact.   
The {\em Tits boundary} of $X$ is denoted by $\tits X$, it is the ideal boundary equipped with 
the Tits metric $|\cdot,\cdot|_T$ -- the intrinsic metric associated to the {\em Tits angle}. 
 The Tits boundary of a CAT(0) space is a CAT(1) space.
If $X$ is a geodesically complete CAT(0) space, then any join decomposition of $\tits X$ is induced by a metric product decomposition of $X$ \cite[Proposition~2.3.7]{KleinerLeeb}.

A {\em $n$-flat} $F$ in a CAT(0) space $X$ is a closed convex subset isometric to $\R^n$.
Its Tits boundary $\tits F\subset\tits X$ is a round $(n-1)$-sphere.
On the other hand, if $X$ is locally compact and $\si\subset \tits X$ is a round  $(n-1)$-sphere, then either there exists an
$n$-flat $F\subset X$ with $\tits F=\si$, or there exists a round $n$-hemisphere $\tau^+\subset\tits X$
with $\si=\D\tau^+$ \cite[Proposition~2.1]{Leeb}. Consequently, if $\tits X$ is $(n-1)$-dimensional, then any round $(n-1)$-sphere in $\tits X$
is the Tits boundary of some $n$-flat in $X$.  A {\em flat ($n$-dimensional) half-space} $H\subset X$ is a closed convex subset isometric to a Euclidean half-space $\R_+^n$.

\subsection{Dynamics at infinity}\label{subsec_dyn}\label{sec_dyn}

Recall that
for any CAT(0) space, the isometry group acts by homeomorphism on the ideal boundary 
and by isometries on the Tits boundary.
Let $X$ be a locally compact CAT(0) space with a {\em geometric group action} $\Ga\acts X$ -- a properly discontinuous cocompact action by isometries.
Then we are interested in the induced action $\Ga\acts\geo X$.

Recall the following construction from \cite{GS_trans}.
Let $G$ be a discrete group acting on a compact Hausdorff space $Z$.
Identify the Stone-Čech compactification $\beta G$ with the set of ultrafilters on 
$G$. For every $\om\in\beta G$ define 
\[T^\om:Z\to Z;\ z\mapsto\wlim_g gz. \] 
Then for fixed $z_0\in Z$ the map $\beta G\to Z$ which sends $\om$
to $T^\om z_0$ is continuous. The family of operators $\{T^\om\}_{\om\in\beta G}$
is closed under composition and the inverse map $g\mapsto g^{-1}$
extends to a continuous involution $S:\beta G\to\beta G$.

Since $X$ is locally compact,
$\bar X=X\cup\geo X$ is a compact Hausdorff space and the above construction
applies. By the semi-continuity of the Tits metric, every operator $T^\om$ is a 1-Lipschitz map $\tits X\to\tits X$ \cite{GS_trans}.

By \cite{Kleiner}, the Tits boundary of $X$ is finite-dimensional. Moreover, $\dim(\tits X)$ coincides with the dimension of a top-dimensional round sphere
$\hat\si\subset\tits X$.
Let $\si\subset\tits X$ be a round sphere and Suppose that for $\om\in\beta\Ga$ exists a rounds sphere $\si\subset\tits X$ such that 
the restriction
$T^\om|_\si$ is an isometry and $T^\om(\tits X)=T^\om(\si)$. Then we say that $\om$ {\em folds} $\si$ onto $T^\om(\si)$ and $T^\om(\si)$ is a {\em folded sphere}. 

\blem[{\cite[Lemma~3.25]{GS_trans}, \cite[Proposition~2.1]{Leeb}}]
For every top-dimensional round sphere $\si\subset\tits X$ there exists $\om\in\beta\Ga$ which folds $\si$.
\elem

A closed subset $M\subset \geo X$ is called {\em minimal}, if it is $\Ga$-invariant and for every $\zeta\in M$
the orbit $\Ga\zeta$ is dense in $M$. If $\geo X$ is not minimal, then by compactness, there exists a non-trivial minimal set $M\subset\geo X$.
Two minimal sets $M, M'$ are either disjoint or coincide.

We say that two points $\xi^\pm\in\geo X$ are {\em $\Ga$-dual} if there exists a sequence $(\ga_k)\subset\Ga$
such that $\ga_k^{\pm 1}o\to\xi^\pm$ for some (and then any) point $o\in X$. The group $\Ga$ satisfies the {\em duality condition}
if the ideal endpoints of any complete geodesic in $X$ are $\Ga$-dual. If $X$ is geodesically complete, then $\Ga$ satisfies the duality condition if and only
if every complete geodesic is nonwandering mod $\Ga$ -- a common condition in dynamics \cite[Corollary~III.1.4]{ballmannbook}.
Recall once again that if a group $G$ acts isometrically, properly discontinuously and with finite covolume on a Hadamard manifold, then $G$
satisfies the duality condition \cite[p.~39]{ballmannbook}. Also, if $\Ga$ acts geometrically on a 2-dimensional geodesically complete CAT(0) polyhedron, then
$\Ga$ satifies the duality condition \cite{BaBr_orbi}. It is expected that any group which acts geometrically on a Hadamard space satisfies the duality condition 
\cite[Question~1.6]{ballmannbook}.

\bprop[{\cite[Proposition~III.1.9]{ballmannbook}}]\label{prop_orbi_min}
If $X$ is a Hadamard space and $\Ga$ is a group of isometries which satisfies the duality condition, then for any $\xi\in\geo X$ the closure of the orbit $\Ga\xi$ is minimal.
\eprop

If $\geo X$ is not minimal, then this leads to a decomposition of $\geo X$ into minimal sets. In particular, if $c\subset X$
is a complete geodesic, then we find two minimal sets $M^\pm\subset\geo X$
such that $c(\pm\infty)\in M^\pm$. Thus, \cite[Lemma~27]{Ri_rank} implies the following.

\blem\label{lem_dual_sym}
If $X$ is a Hadamard space and $\Ga$ satisfies the duality condition, then $\geo X$ is either minimal or contains a non-trivial
closed symmetric set.
\elem

Recall from \cite{GS_trans}:
\bdfn
A pair of points $\xi, \eta\in\geo X$ is {\em compressible} if 
\[|T^\om \xi, T^\om \eta|<|\xi,\eta|\]
for some $\om\in\beta\Ga$.
 A subset $A\subset\geo X$ is {\em compressible} if it contains a compressible pair.
Otherwise $A$ is called {\em incompressible}.
\edfn

Note that for every $\om\in\beta\Ga$ the operator $T^\om$ restricts to an isometry on every incompressible set.

\subsection{Morse quasiflats and Morse rank}\label{sec_morse}

Morse quasiflats are a higher dimensional generalization of Morse quasi-geodesics. They were originally introduced and studied in \cite{HKS_I,HKS_II}. 
They are quasi-isometry invariant and  display nice asymptotic behavior, namely they 
have unique tangent cones at infinity. Here we use them to define an algebraic invariant of a discrete group.
We begin with a very particular case of Morse quasiflats, namely periodic Morse flats. Recall that a flat $F$ in a metric space $X$ is {\em periodic}, if its stabilizer in the full isometry group of $X$ contains a subgroup which acts geometrically on $F$.
A periodic $n$-flat $F$ is called {\em Morse}, if it does not bound a flat $(n+1)$-half-space.
Periodic Morse $1$-flats are precisely Ballmann's rank 1 axes \cite[Chapter~III.3]{ballmannbook}. If  
$\ga$ is a periodic geodesic in a CAT(0) space $X$, then $\ga$ is Morse if and only if (one and then) both ideal endpoints of
$\ga$ are isolated in the Tits boundary. Similar characterizations hold for Morse quasiflats 
(\cite[Corollary~11.5]{HKS_I},\cite[Proposition~6.14]{HKS_II}).
Periodic Morse flats are important because they lead to many regular points in the Tits boundary.

\bthm[{\cite[Theorem~D]{St_rrI}}]\label{thm_reg}
Let $X$ be a locally compact CAT(0) space.  Suppose that $X$ contains a periodic Morse flat $F$. Then 
the complement in $\tits F$ of the set of regular points can be covered by a finite set of round spheres of positive codimension.
\ethm

\bdfn
Let $\Ga$ be a finitely generated discrete group.
A free abelian subgroup   $G<\Ga$ of rank $n$ is called {\em Morse}, if it is a quasiflat and for every $L\geq 1$, $A\geq 0$, $r\geq 0$
there exists $R> 0$ such that any $n$-dimensional $(L,A)$-quasidisc $D\subset\Ga$ with boundary $\D D$ in the $r$-neighborhood of $G$
lies entirely in the $R$-neighborhood of $G$. The {\em Morse rank} of $\Ga$ is the smallest natural number $k$ such that $\Ga$
contains a Morse free abelian subgroup of rank $k$. If no such subgroup exists then the Morse rank is infinite. We denote the Morse rank of $\Ga$
by $\rank_{Mo}(\Ga)$. 
\edfn

Note that if $M$ is a Hadamard manifold with a geometric group action $\Ga\acts M$, then Rank Rigidity implies that 
\[\rank(M)=\rank_{Mo}(\Ga)\]
where the rank of $M$ is defined as the largest natural number $k\in \N$ such that every geodesic in $M$ admits $k$ linearly independent parallel Jacobi fields.

The following simple observation makes the Morse rank useful.

\bprop\label{prop_morse_rank}
Let $X$ be a Hadamard space with a geometric action $\Ga\acts X$. If the Morse rank of $\Ga$ is finite,  $\rank_{Mo}(\Ga)=n\in\N$, then
$X$ contains a $\Ga$-periodic Morse $n$-flat. In particular, the Tits boundary of $X$ contains regular points.
\eprop

\proof
Let $G<\Ga$ be a Morse free abelian subgroup of rank $n$. By the Flat Torus Theorem,  $X$ contains a $G$-periodic $n$-flat $F$.
Since $G$ is Morse  and $\Ga$ is quasi-isometric to $X$, we see that $F$ is Morse as well. The supplement follows from Theorem~\ref{thm_reg}.
\qed

\section{A submetry at infinity}\label{sec_sub}

\subsection{Symmetric subsets}\label{sec_sym}

\bdfn
Let $X$ be a  CAT(0) space.
Two points $\xi,\hat\xi\in\geo X$ are called {\em antipodal}, if $|\xi,\hat\xi|\geq\pi$. They are called
{\em visually antipodal}, if there is a complete geodesic $c$ in $X$ with $\geo c=\{\xi,\hat\xi\}$. 
For a subset $A\subset\geo X$ we define the set of (visually) antipodal points of $A$ by $\Ant_{(vis)}(A)$.
We say $A$ is {\em (visually) symmetric}, if 
\[\Ant_{(vis)}(A)\subset A.\]
\edfn

Visually symmetric subsets are also called {\em involutive} \cite{E_sym}.
Recall that  the closure of a visually symmetric subset $A\subset\geo X$ is symmetric \cite[Corollary~3.6]{St_rrII}.

Recall that 
a point $\eta$ in a CAT(1) space $Z$ is {\em almost near} to a point $\xi\in Z$, if for every antipode $\hat\xi$ of $\xi$ holds
\[|\xi,\eta|+|\eta,\hat\xi|=\pi.\] 
No geodesic between almost near points can contain a branch point in the interior \cite[Lemma~3.16]{St_rrII}.

\subsection{Minimal closed symmetric sets}

The following result is a consequence of \cite[Proposition~4.4, Proposition~4.6]{St_rrII}.
\bthm\label{thm_submetry}
Let $X$ be a locally compact and geodesically complete CAT(0) space.
If $\geo X$ contains a non-trivial closed symmetric subset, then the
family of minimal closed symmetric subsets in $\geo X$ constitute an   equidistant decomposition of $\tits X$ and therefore induce
a submetry 
\[\delta:\tits X\rightarrow \Delta.\]
The following additional properties hold.
\begin{enumerate}
	\item $\tits X$ has diameter $\pi$ and is horizontally geodesically complete. More precisely, any pair of $\delta$-near points is almost near and every horizontal geodesic lies in a round circle.
	\item $\delta:\geo X\to\Delta$ is $\Isom(X)$-equivariant.
	\item $\delta$ is transnormal, i.e. every extension of a horizontal geodesic is piecewise horizontal.
\end{enumerate}
\ethm

For the rest of this section we keep the setting and notation of Theorem~\ref{thm_submetry}.

\blem[{\cite[Lemma~5.14]{St_rrII}}]\label{lem_restrsub}
Let $\si\subset\tits X$ be a round sphere of dimension $d$. Suppose that every for every point $\xi\in\si$ and every direction $v\in\Si_\xi\tits X$ 
there exists an antipode $\hat v$ in 
$\Si_\xi\si$. Then the restriction $\delta_\si:=\delta|_\si$ is a submetry $\delta_\si:\si\to\Delta$ and every $\delta_\si$-horizontal geodesic is
also $\delta$-horizontal. In particular, $\Delta$ is a compact Alexandrov space of dimension at most $d$ and curvature at least $1$. 
\elem
\medskip
Note that the above lemma applies to all  top-dimensional round spheres $\si\subset\tits X$ \cite[Lemma~2.1]{BL_building}.

\blem\label{lem_cont} 
The map $\delta:\geo X\to\Delta$ is continuous.
\elem

\proof
Let $\xi_k\to\xi$ be a converging sequence in $\geo X$ and set $x_k:=\delta(\xi_k)$.
Since $\Delta$ is compact, we can pass to a convergent subsequence $x_{k_l}\to x$.
Denote by $B_k\subset\geo X$ the minimal closed symmetric subset containing $\xi_k$, i.e. the $\delta$-fiber of $x_k$.
Let $B$ be the fiber of $x$.
By definition, $B_{k_l}$ converges to $B$ with respect to the Hausdorff metric induced by the Tits metric.
By semi-continuity of the Tits metric, it follows that $\xi\in B$, i.e. $\delta(\xi)=x$.
Thus $\delta$ is continuous.
\qed

\medskip

For a point $\xi\in\tits X$ we define the {\em horizontal space} $H_\xi\subset\Si_\xi\tits X$
as the subset of starting directions of horizontal geodesics.

Recall that if $Z$ is a CAT(1) space with a point $z\in Z$, then we call a subspace $H\subset\Si_z Z$ {\em almost symmetric},
if for every point $v\in H$ all genuine antipodes $\hat v\in \Si_z Z$ are again contained in $H$.
From \cite[Proposition~4.12, Lemma~4.15, Lemma~4.16]{St_rrII} we have:

\blem\label{lem_cloconsym}
Every horizontal space $H_\xi\subset\Si_\xi\tits X$ is closed, convex and almost symmetric.
Moreover, for every $\xi\in\tits X$ there exists $\eps>0$ such that for each $v\in H_\xi$ there is a unique geodesic of length $\eps$ in the direction $v$, and this geodesic
 is $\delta$-horizontal.
\elem

\subsection{Singular spheres}

\bdfn
Let $\hat\si\subset\tits X$ be a top-dimensional round sphere.
A round sphere $\si\subset\hat\si$ is called {\em singular}, if there exists a round hemisphere $\tau^+\subset\tits X$
with $\D\tau^+=\si$ and  $\tau^+\cap\hat\si=\si$ such that the union $\tau^+\cup\tau^-$ is a round sphere for every 
round hemisphere $\tau^-\subset\hat\si$ with $\D\tau^-=\si$.
\edfn

\blem\label{lem_refl_sph}
Let $\hat\si\subset\tits X$ be top-dimensional round sphere and $\si\subset\hat\si$ a singular sphere.
Then $\si$ is a reflecting sphere for the submetry  $\delta_{\hat\si}$.
\elem

\proof
Let $\tau\subset\tits X$ be a round hemisphere for $\si$ as in the definition of singular spheres.
For every point $\xi\in\hat\si$ has an antipode $\xi_1\in\tau$. Then every hemisphere $\tau'\subset\hat\si$
with $\D\tau'=\si$ contains an antipode $\xi_2$ of $\xi_1$. This shows the claim since the fibers of $\delta_{\hat\si}$ are the intersections of 
minimal closed symmetric sets with $\hat\si$ . 
\qed

\section{Symmetries of induced submetries}\label{sec_ind_sub}

For this section we will assume that $X$ is a cocompact Hadamard space and that $\geo X$
contains a non-trivial closed symmetric subset. We will denote by $\delta:\tits X\to\Delta$
the induced submetry provided by Theorem~\ref{thm_submetry}.

\subsection{Producing singular spheres}

\blem\label{lem_orth_half}
Let  $F\subset X$ be a flat with $\si:=\tits F$ and $\xi\in\si$. 
Suppose there exists a sequence  $(\xi_k)\subset\tits X\setminus\si$ with $\xi_k\to\xi$ in $\tits X$.
Then there exists a flat $F'\subset X$ with $\si':=\tits F'$ and 
a sequence $(\ga_k)$ in $\Ga$ such that $\ga_k|_\si$ converges pointwise to an isometry $\psi:\si\to\si'$. 
Moreover, 
there is  a flat half-plane $H'\subset X$ orthogonal to $F'$ such that
$\psi(\xi)\in\tits H\cap\si'$.
\elem

\proof
Choose a base point $p\in F$.
For each $k\in\N$ choose $x_k\in p\xi_k$ such that $|x_k,F|=k$. Denote by $\bar x_k$ the nearest-point projection of $x_k$ to $F$.
Choose a sequence $(\ga_k)\subset\Ga$ such that after possibly passing to subsequences we have the following convergences: 
$\ga_k(F,\bar x_k)\to (F',p')$, $\ga_k|_\si$ converges pointwise to an isometry $\psi:\si\to\si':=\tits F'$, 
$\ga_k x_k\to\eta\in\geo X$,
$\ga_k\xi_k\to\zeta\in\si'$ and $\ga_k p\to\hat\zeta\in\si'$. We claim that $|\zeta,\eta|=|\eta,\hat\zeta|=\frac{\pi}{2}$.
Suppose for the moment that this is true. Then the claim follows since $\zeta$ and $\hat\zeta$ are antipodes and the geodesic ray $p'\eta$
is orthogonal to $F'$. To prove the claim let $y\in p'\eta$ and choose a sequence $(y_k)$ with $y_k\in x_k\bar x_k$ such that 
$\ga_k y_k\to y$.
Note that $\angle_{\bar x_k}(\xi,x_k)\geq\frac{\pi}{2}$ and $\angle_{\bar x_k}(p,x_k)\geq\frac{\pi}{2}$. Thus 
$\angle_{x_k}(\bar x_k,p)\leq\frac{\pi}{2}$  and $\angle_{x_k}(\bar x_k,\xi_k)\geq\frac{\pi}{2}$.
Therefore  
\[\angle_{y_k}(x_k,\xi_k)\leq\frac{\pi}{2},\ \angle_{y_k}(\bar x_k,\xi_k)\geq\frac{\pi}{2},\ \angle_{y_k}(\bar x_k,p)\leq\frac{\pi}{2}
,\ \angle_{y_k}(x_k,p)\geq\frac{\pi}{2}.\]
Taking limits and using semi-continuity of angles, we obtain $\angle_{p'}(\zeta,y)=\angle_{y}(\zeta,p')=\frac{\pi}{2}$.
In particular,
\[\lim_{k\to\infty}\angle_{y_k}(\bar x_k,\xi_k)=\lim_{k\to\infty}\angle_{\bar x_k}(y_k,\xi)=\frac{\pi}{2}.\]
Using the estimates above and semi-continuity of angles again, we deduce
\[\lim_{k\to\infty}\angle_{x_k}(\bar x_k,\xi_k)=\lim_{k\to\infty}\angle_{x_k}(\bar x_k,p)=\frac{\pi}{2}\]
and
\[\lim_{k\to\infty}\angle_{y_k}(x_k,p)=\lim_{k\to\infty}\angle_{y_k}(\bar x_k,p)=\frac{\pi}{2}.\]
Thus,  $\angle_{p'}(\hat\zeta,y)=\angle_{y}(\hat\zeta,p')=\frac{\pi}{2}$. In particular, the infinite triangles $\Delta(\zeta,p',y)$ and $\Delta(\hat\zeta,p',y)$
are flat strips. Since this holds for all $y\in p'\eta$, the claim follows.
\qed
\medskip

Recall that a {\em spherical $m$-lune} in a CAT(1) space is a closed convex subset $\la$ isometric to 
$S^{m-1}\circ[0,\theta]$ for some $\theta\in(0,\pi]$.

\blem\label{lem_orth_lune}
Let $\tilde\si\subset\tits X$ be a top-dimensional round sphere and $\tilde\la\subset\tits X$ a spherical $m$-lune spanned by round hemispheres $\tilde\tau^\pm\subset\tits X$
such that $\tilde\la\cap\tilde\si=\tilde\tau^-$.
Then there exists a top-dimensional round sphere $\hat\si\subset\tits X$ 
and a sequence $(\ga_k)$ in $\Ga$ such that $\ga_k|_{\tilde\si}$ converges pointwise to an isometry $\psi:\tilde\si\to\hat\si$ with the following properties.
There is a round  $m$-hemisphere $\hat\tau\subset\tits X$ with $\hat\tau\cap\hat\si=\D\hat\tau$ and such that 
 $\psi(\tilde\tau^-)\subset\D\hat\tau$.
\elem

\proof
Since $\tilde\si\subset\tits X$ is top-dimensional, there exists a flat $\tilde F\subset X$ with $\tits \tilde F=\tilde\si$.
We have $\tilde\la=\D\tilde\tau^-\circ [0,\theta]$. Write $\tilde F=\tilde F_1\times \tilde F_2$ with $\tits \tilde F_1=\D\tilde\tau^-$.
Then the flat $\tilde F_2$ lies in the cross section $CS(\tilde F_1)$ of the parallel set of $\tilde F_1$.
Therefore, the Tits boundary $\tits CS(\tilde F_1)$ contains $\tits\tilde F_2$ and the arc $[0,\theta]$.
Since $\tits\tilde F_2$ and $[0,\theta]$ intersect only in one point, the claim follows similarly as in Lemma~\ref{lem_orth_half}.
Just observe that if $p\in F$ is point and $\xi\in\tits\tilde\la$ is a point in the cross section $[0,\theta]$ of $\tilde\la$, 
then the ray $p\xi$ lies in the parallel set $P(F_1)$. We can then produce a limit flat $F'=F'_1\times F_2'$ and a flat half-plane
$H'\subset X$ which is orthogonal to $F'$ and contained in the parallels set $P(F_1')$. Then $\hat\tau:=\tits (H'\times F_1')$
is the required  round  $m$-hemisphere.
\qed

\bprop\label{prop_sing_sph}
Let $\si\subset\tits X$ be a top-dimensional round sphere and $\xi\in\si$. 
Suppose there exists a sequence  $(\xi_k)\subset\tits X\setminus\si$ with $\xi_k\to\xi$ in $\tits X$.
Then there exists a top-dimensional round sphere $\hat\si\subset\tits X$ and a sequence $(\ga_k)$ in $\Ga$
such that $\ga_k|_\si$ converges pointwise to an isometry $\psi:\si\to\hat\si$. Moreover, $\hat\si$ contains a singular sphere $\si$ 
and $\psi(\xi)\in\si$.
\eprop

\proof
We apply Lemma~\ref{lem_orth_half} to find a top-dimensional round sphere $\si_1\subset\tits X$ and an isometry $\psi_1:\si\to\si_1$ in the closure of 
$\Ga<\Isom(X)$ and a round 1-hemisphere $\tau_1$ with $\tau_1\cap\si_1=\D\tau_1$ and $\psi_1(\xi)\in\D\tau_1$. 
By the Lune Lemma (Lemma~\ref{lem_lune}), $\D\tau_1$ is either singular in $\si_1$ and we are done, or it spans a spherical 2-lune $\la_1$ with some round 1-hemisphere $\tau_1^-\subset\si_1$.
In the latter case we apply Lemma~\ref{lem_orth_lune} to find a top-dimensional round sphere $\si_2\subset\tits X$ and an isometry $\psi_2:\si_1\to\si_2$ in the closure of 
$\Ga<\Isom(X)$ and a round 2-hemisphere $\tau_2$ with $\tau_2\cap\si_2=\D\tau_2$ and $\psi_2(\tau_1^-)\subset\D\tau_2$. 
Again, by the Lune Lemma, $\D\tau_2$ is either singular in $\si_2$ and we are done, or it spans a spherical 3-lune $\la_2$ with some round 2-hemisphere $\tau_2^-\subset\si_2$. Proceeding in this manner, we arrive at the conclusion in at most $d$ steps.
\qed

\subsection{Symmetry vs. regularity}

\blem\label{lem_limit_sym}
Let $\si,\si'\subset\tits X$ be top-dimensional round spheres.
Let $\varphi$ be an isometry of $\si'$ such that $\delta_{\si'}\circ\varphi=\delta_{\si'}$. 
Suppose that $(\ga_k)\subset \Ga$ is a sequence, such that $\ga_k|_\si$ converges pointwise to an isometry $\psi:\si\to\si'$.
Then  $\delta_\si=\delta_\si\circ\psi^{-1}\circ\varphi\circ\psi$.
\elem

\proof
After passing to a subsequence we may assume that the isometries 
$\ga_k:\Delta\to\Delta$ converge uniformly to an isometry $g:\Delta\to\Delta$. Using continuity and $\Isom(X)$-equivariance of $\delta$ we obtain
for every $\xi\in\si$:
\begin{align*}
&g\delta(\xi)=\lim_{k\to\infty}\ga_k(\delta(\xi))=\delta(\lim_{k\to\infty}\ga_k\xi)=\delta(\psi(\xi))=\\
&\delta(\psi(\psi^{-1}\circ\varphi\circ\psi(\xi)))=\delta(\lim_{k\to\infty}\ga_k(\psi^{-1}\circ\varphi\circ\psi(\xi))=\\
&\lim_{k\to\infty}\ga_k\delta(\psi^{-1}\circ\varphi\circ\psi(\xi))=g\delta(\psi^{-1}\circ\varphi\circ\psi(\xi)).
\end{align*}
Thus $\delta_\si=\delta_\si\circ\psi^{-1}\circ\varphi\circ\psi$ as required.  
\qed

\medskip
Recall that a point $\xi$ in a top-dimensional round sphere $\si\subset\tits X$ is {\em singular}
if and only if there exists a sequence  $(\xi_k)\subset\tits X\setminus\si$ with $\xi_k\to\xi$ in $\tits X$.
Otherwise $\xi$ is a regular point since it has a neighborhood in $\tits X$ which lies entirely in $\si$.

\bprop\label{prop_sing_refl}
Let $\hat\si\subset\tits X$ be a top-dimensional round sphere.
Then for every singular point $\xi\in\hat\si$ there exists a reflecting sphere $\si\subset\hat\si$ for $\delta$
with $\xi\in\si$.
\eprop

\proof
By Proposition~\ref{prop_sing_sph}, there exists a top-dimensional round sphere $\hat\si'\subset\tits X$
and an isometry $\psi:\hat\si\to\hat\si'$ which is a pointwise limit of a sequence $(\ga_k|_\si)$ with $(\ga_k)\subset\Ga$, such that $\hat\si'$
contains a singular sphere $\si'$ and $\psi(\xi)\in\si'$.
By Lemma~\ref{lem_refl_sph}, $\si'$ is a reflecting sphere for $\delta|_{\hat\si'}$.
By Lemma~\ref{lem_limit_sym}, $\si:=\psi^{-1}(\si)$ is a reflecting sphere for $\delta|_{\hat\si}$.
\qed

\section{Irreducible spaces and geodesically complete Tits boundaries }\label{sec_irr}

We keep working under the assumptions of the last section: $X$ is a cocompact Hadamard space with a non-trivial
closed symmetric subset in $\geo X$ and $\delta:\tits X\to\Delta$ is the induced submetry.
In this section we aim to show that if $X$ is irreducible, then the dimension of the base $\Delta$ has to coincide with the dimension of the 
Tits boundary $\tits X$.

\subsection{Splittings via non-branching piecewise horizontal geodesics }

\bprop\label{prop_split}
Let $\hat\si\subset\tits X$ be a folded round sphere. Suppose  that there exists a 
round sphere $\si\subset\hat\si$ which splits as a 
spherical join
$\si=\si^-\circ\si^+$ such that every geodesic $\xi^-\xi^+$ with $\xi^\pm\in\si^\pm$ is piecewise horizontal and 
does not branch in $\tits X$. Further, assume that every normal horizontal direction at $\si^-$ points to $\si^+$.
Then $X$ splits as a non-trivial product.
\eprop

\proof
Set $A:=\Ant(\si^-)\subset\geo X$. We will show that $A$ is closed,  symmetric and convex.

We first claim that a point $\xi\in \tits X$ lies in $A$ if and only if $\xi$ is the center of a round hemisphere $\tau\subset\tits X$ with $\D\tau=\si^-$
and such that the geodesic $\xi\xi^-$ is piecewise horizontal for every $\xi^-\in\si^-$.

It is enough to show that a point $\xi\in\tits X$ with this property
has to lie in $A$. Let $\eta\in\si^-$ be an arbitrary point and extend the
piecewise horizontal geodesic 
$\xi\eta$ to an antipode $\hat\xi\in\hat\si$. Then we have  $|\hat\xi,\zeta|=\frac{\pi}{2}$ for every $\zeta\in\si^-$.
Thus, the piecewise horizontal geodesic $\eta\hat\xi\subset\hat\si$ is normal to $\si^-$. By assumption, $\hat\xi$ has to lie in $\si^+$ and therefore 
$\xi\in A$ as claimed. 

\noindent\textit{Claim 1:} $A\subset\geo X$ is closed.

Let $(\xi_k)\subset A$ with $\xi_k\to\xi$.
By semi-continuity of the Tits metric, we have $|\xi,\eta|=\frac{\pi}{2}$ for every $\eta\in\si^-$.
The lune lemma implies that
$\xi$ is the center of a round hemisphere $\tau\subset\tits X$ with $\D\tau=\si^-$.
For every  $\eta\in\si^-$  the directions $v_k\in\Si_\eta\tits X$ of $\eta\xi_k$ are horizontal and converge to the direction
$v\in\Si_\eta\tits X$ of $\eta\xi$.
It follows from Lemma~\ref{lem_cloconsym} that $v$ is horizontal. Thus $\xi$ lies in $A$.

\noindent\textit{Claim 2:} $A\subset\geo X$ is symmetric.

Let $\xi\in A$ and $\hat\xi\in\Ant(\xi)$. Note that since geodesics $\xi^-\xi^+\subset\si$ do not branch in $\tits X$ for $\xi^\pm\in\si^\pm$,
none of the geodesics $\xi\xi^-$ with $\xi^-\in\si^-$ can branch in $\tits X$. Thus, for every $\xi^-\in\si^-$ the piecewise horizontal geodesic 
$\xi\xi^-$ extends to a geodesic to $\hat\xi$. By the lune lemma, $\hat\xi$ is the center of a round hemisphere $\hat\tau\subset\tits X$ with $\D\hat\tau=\si^-$.
For every $\xi^-\in\si^-$  the geodesic $\hat\xi\xi^-$ is piecewise horizontal and therefore $\hat\xi\in A$.

\noindent\textit{Claim 3:} $A\subset\tits X$ is convex.

Let $\xi$ and $\xi'$ be points in $A$ at distance  $\theta<\pi$.
By the lune lemma, the geodesic $\xi\xi'$ and $\si^-$  span a spherical lune $\la\cong\si^-\circ[0,\theta]$ in $\tits X$.
Let $\zeta\in\la$ be the midpoint of $\xi$ and $\xi'$. 
For $\eta\in\si^-$ denote by $v, v'\in\Si_\eta\tits X$ the directions of the piecewise horizontal geodesic $\eta\xi$ and $\eta\xi'$ respectively.
Note that $|v,v'|=\theta$.
By Lemma~\ref{lem_cloconsym}, the geodesic $\eta\zeta$ is piecewise horizontal. Since this holds for every point $\eta\in\si^-$, we obtain $\zeta\in A$.

\noindent\textit{Claim 4:} $X$ splits isometrically.

Since $A\subset\geo X$ is closed, convex and symmetric, \cite[Theorem~2.9]{BMS_affine} implies that the Tits boundary has to split as a spherical join, 
$\tits X\cong A\circ A^\perp$. Thus the claim follows from \cite[Proposition~2.3.7]{KleinerLeeb}.
\qed


\bcor\label{cor_point_split}
Let $\hat\si\subset\tits X$ be a folded round sphere. Suppose  that there exists a point $\zeta\in\hat\si$
such that piecewise horizontal geodesics in $\hat\si$ which start in $\zeta$ do not branch in $\tits X$
before time $\frac{\pi}{2}$. 
Then $X$ splits as a non-trivial product.
\ecor

\proof
Let $\hat\zeta\in\hat\si$ be the antipode of $\zeta$. Set $\si^-:=\{\zeta,\hat\zeta\}$
and $\si^+:=\exp_{\frac{\pi}{2}}(H_{\zeta}\cap\Si_\zeta\hat\si)$. By Lemma~\ref{lem_cloconsym}, $\si^+\subset\hat\si$
is a round sphere and Proposition~\ref{prop_split} applies.
\qed

\bcor\label{cor_codim1_split}
Let $\hat\si\subset\tits X$ be a folded round sphere. Suppose $\si\subset\hat\si$
is a reflecting sphere of codimension $2$. Then $X$ splits as a non-trivial product.
\ecor

\proof
Write $\hat\si=\si\circ\si^\perp$ and let $\zeta^+\in\si^\perp$.
Since $\delta$ is constant on $\si^\perp$, the $\delta_{\hat\si}$-horizontal space $H_{\zeta^+}\cap\Si_{\zeta^+}\hat\si$
has to be orthogonal to $\si^\perp$. 
We claim that piecewise horizontal geodesics in $\hat\si$ which start in $\zeta^+$ cannot branch in $\tits X$
before time $\frac{\pi}{2}$.

To see this, let $\zeta^+\xi\subset\si$ be such a geodesic and let $\hat\xi\in\hat\si$ be the antipode of $\xi$. 
Denote by $c\subset\hat\si$  the 
geodesic from $\xi$
to $\hat\xi$ through $\zeta^+$. Let $0<t_1<t_2<\ldots<t_l=\pi$ be minimal such that $c|_{[t_{i-1},t_{i}]}$
is horizontal, $1\leq i\leq l$. Rotating $c$  around $\si$
we obtain a round $2$-sphere $\tilde\si\subset\hat\si$ such that every geodesic $\tilde c\subset\tilde\si$
from $\xi$
to $\hat\xi$ is horizontal on $[t_{i-1},t_{i}]$
for $1\leq i\leq l$. It follows from \cite[Lemma~4.20]{St_rrII} that each $\tilde c$
can only branch at time $t=\frac{\pi}{2}$.
Thus the claim follows from Corollary~\ref{cor_point_split}.
\qed

\bcor\label{cor_inf_split}
Let $\hat\si\subset\tits X$ be a folded round sphere.
Let $\{\si_\al\}_{\al\in\mathcal A}$ be the family of reflecting spheres in $\hat\si$.
If $\mathcal A$ is infinite, then $X$ splits non-trivially as a product.
\ecor

\proof
If $\mathcal A$ is infinite, then we find a non-trivially converging sequence $\si_k\to\si$ of reflecting hyperspheres.
We may assume that the codimension 2 spheres $\si_k\cap\si$ also converge to a codimension 2 sphere $\tilde\si\subset\hat\si$.
Since the angle between $\si_k$ and $\si$ goes to zero, and $\delta$ is invariant under reflection at $\si_k$ and $\si$, 
we see that $\tilde\si\subset\hat\si$ is a reflecting sphere of codimension 2. The claim follows from Corollary~\ref{cor_codim1_split}.
\qed

\blem\label{lem_euc_fac}
Let $\hat\si\subset\tits X$ be a top-dimensional round sphere.
Let $\{\si_\al\}_{\al\in\mathcal A}$ be the family of reflecting spheres in $\hat\si$.
Suppose $\mathcal A$ is finite and $\bigcap_{\al\in\mathcal A}\si_\al$ is non-empty. Then $X$ splits non-trivially as a product.
\elem

\proof
By Proposition~\ref{prop_sing_refl} we know that any point $\xi\in\hat\si$ which does not lie on a reflecting sphere has to be regular.
If there exists a piecewise horizontal geodesic in $\hat\si$ which runs through a regular point and intersects two reflecting hyperspheres in angles different from $\frac{\pi}{2}$,
then this geodesic spans a $2$-sphere $\si_H\subset\hat\si$ all of whose geodesics are piecewise horizontal. Since $\si_H$ has a regular point, 
it provides a regular piecewise horizontal geodesic  of length $\pi$. The endpoints $\xi,\hat\xi$ of this geodesic are each others only antipodes
in $\tits X$.
As $X$ is geodesically complete, it has to coincide with the parallel set $P(\xi,\hat\xi)$, in particular, $X$ has to split off a line.
 
Hence, we may assume that all reflecting hyperspheres intersect orthogonally.
In this case, we find a regular piecewise horizontal geodesic  of length $\frac{\pi}{2}$ (through every regular point).
Let $\zeta^+\xi\subset\hat\si$ be such a geodesic. We claim that piecewise horizontal geodesics in $\hat\si$
which start in $\zeta^+$ do not branch in $\tits X$ before time $\frac{\pi}{2}$. Note that in our setting, the directions in $H_{\zeta^+}\cap\Si_{\zeta^+}\hat\si$
which point to geodesics that run in the regular part of $\hat\si$ up to time $\frac{\pi}{2}$ form an open and dense set.
Since these geodesics cannot branch before time $\frac{\pi}{2}$, horizontal geodesic completeness of $\tits X$ and a limit argument show the claim.
Thus, the statement follows from Proposition~\ref{prop_split}.
\qed

\subsection{Irreducibility and the dimension of the base}

\blem\label{lem_fin_fac}
Let $\hat\si\subset\tits X$ be a folded round sphere.
Let $\{\si_\al\}_{\al\in\mathcal A}$ be the family of reflecting spheres in $\hat\si$.
Suppose $\mathcal A$ is finite and $\bigcap_{\al\in\mathcal A}\si_\al=\emptyset$. Then $X$ splits non-trivially as a product or $\dim(\hat\si)=\dim(\Delta)$.
\elem

\proof
Suppose $\dim(\Delta)<\dim(\hat\si)$.
We infer from Proposition~\ref{prop_splitsub} that there exists a disjoint decomposition $\mathcal A=\mathcal A^+\cup \mathcal A^-$ and a splitting 
$\hat\si=\si^-\circ\si^+$ with $\si^\pm=\bigcap_{\al\in\mathcal A^\pm}\si_\al$ and such that every geodesic $\xi^-\xi^+$
is piecewise horizontal, $\xi^\pm\in\si^\pm$. 
Note that a reflecting sphere $\si_{\al^\pm}$, $\al^\pm\in \mathcal A^\pm$, intersects $\si^\mp$ in codimension 1.
Thus every geodesic in $\hat\si$ between $\si^-$ and $\si^+$ can be approximated by  regular such geodesics.
It follows that none of the piecewise horizontal geodesics $\xi^-\xi^+$, $\xi^\pm\in\si^\pm$, branches in $\tits X$.
By Proposition~\ref{prop_split}, $X$ has to split as a non-trivial product.
\qed

\medskip

From Corollary~\ref{cor_inf_split}, Lemma~\ref{lem_euc_fac}, Lemma~\ref{lem_fin_fac} and Lemma~\ref{lem_polytope} we obtain:   

\bprop\label{prop_sph_orbi}
If $X$ is irreducible, then $\dim(\tits X)=\dim(\Delta)$ and $\Delta$ is a spherical orbifold.
\eprop

\bprop\label{prop_geo_complete}
If $X$ is irreducible, then $\tits X$ is geodesically complete.
\eprop

\proof
By \cite[Lemma~2.4]{St_rrI}, it is enough to show that every point in $\tits X$ lies in a top-dimensional round sphere.
By Proposition~\ref{prop_sph_orbi}, we have $d:=\dim(\tits X)=\dim(\Delta)$ and $\Delta$ is a spherical $d$-orbifold.
Let $\hat\si\subset\tits X$ be a round $d$-sphere. Let $\xi\in\hat\si$ be a point with $\delta(\xi)\in\Delta_{reg}$.
Let $\hat\xi\in\tits X$ be an antipode of $\xi$. If $U\subset\hat\si$ is a neighborhood of $\xi$ such that $\delta$
restricts to an isometric embedding on $U$, then $\xi$ is the only antipode of $\hat\xi$ in $U$. By \cite[Lemma~2.11]{CM_bieber},
there exists a round $d$-sphere $\tilde\si\subset\tits X$ containing $\{\xi,\hat\xi\}$.
Since every point $\zeta\in\tits X$ has an antipode in $\hat\xi$ (\cite[Lemma~2.1]{BL_building}), we see that every point $\zeta\in\tits X$
with $\delta(\zeta)\in\Delta_{reg}$ lies in a round $d$-sphere.
For an arbitrary point $\xi\in\tits X$ choose a sequence $(\xi_k)\subset\geo X$ with $\delta(\xi_k)\in\Delta_{reg}$, $k\in\N$, and $\xi_k\to\xi$.
 After passing to a subsequence, we may assume that 
$\sum_{k=1}^\infty |x_k,x_{k+1}|<\infty$ where $x_k:=\delta(\xi_k)$, $k\in\N$.  
Denote by $c_k\subset\Delta$ a piecewise geodesic made of segments $x_{i-1}x_i$, $i=2,\ldots,k$.
Further, let $c=\bigcup_{k\in\N} c_k$. For every $k\in\N$ we lift $c_k$ to a piecewise horizontal geodesic $\hat c_k$
starting at the point $\xi_k$. Passing to a subsequence, we obtain a horizontal lift $\hat c$ of $c$ as a limit of the curves $\hat c_k$.
Note that $\hat c$ contains $\xi$.
By the first part, we find a round $d$-sphere $\tilde \si\subset\tits X$ which contains a point of $c$, different from $\xi$.
Since $c_k\subset\Delta_{reg}$ 
for every $k\in\N$, we conclude $\hat c\subset\tilde\si$. In particular, $\xi\in\tilde\si$.
Hence the claim follows.
\qed

\section{Proof of main result and geometric applications}\label{sec_main}

\proof[Proof of Theorem~\ref{thm_mainA}]
If $\geo X$ contains a non-trivial closed symmetric subset, then the minimal closed symmetric subsets provide an equidistant decomposition of $\tits X$
and therefore a submetry
\[\delta:\tits X\to \Delta.\]
By Proposition~\ref{prop_geo_complete}, if $X$ is irreducible, then the Tits boundary  $\tits X$ has to be geodesically complete.
Thus, by \cite[Main Theorem]{Ly_rigidity}, it has to be an irreducible spherical building. Consequently, $X$ has to be a Riemannian symmetric space of higher rank or
a Euclidean building of higher dimension \cite[Main Theorem]{Leeb}.
\qed

\proof[Proof of Theorem~\ref{thm_mainB}]
If $\Ga$ acts minimally on $\geo X$, then the geodesic flow of $X$ is topologically transitive mod $\Ga$ \cite[Theorem~III.2.3]{ballmannbook}.
Thus it has to have a dense orbit \cite[Proposition~3.4]{E_conj}. 
Suppose $\geo X$ contains a non-trivial closed $\Ga$-invariant subset. Then by Lemma~\ref{lem_dual_sym}, $\geo X$ contains a non-trivial
closed symmetric set. Thus the claim follows from Theorem~\ref{thm_mainA}.
\qed
\medskip

Before we turn to Theorem~\ref{thm_mainC}, we need a criterion to ensure non-minimality of an action $\Ga\acts\geo X$.

\bprop\label{prop_reg_dich}
Let $X$ be a Hadamard space with a geometric group action $\Ga\acts X$. Suppose that the Tits boundary
contains a relative compact open subset $U\subset\tits X$.
Then either $X$ contains a $\Ga$-periodic rank 1 axis, or the action $\Ga\acts\geo X$ is not minimal.
\eprop

\proof
Suppose that $X$ does not contain a $\Ga$-periodic rank 1 axis.
Since $U\subset\tits X$ is relatively compact, \cite[Proposition~2.7]{BaBu_periodic} shows that minimality can only fail if
$\tits X$ contains a pair of distinct points $\zeta,\zeta'$ such that for every point $x\in X$ the sector spanned the geodesic rays $x\zeta$ and $x\zeta'$ is flat.
But then $\{\zeta,\zeta'\}$ is an incompressible pair.
By \cite[Section~5]{GS_trans}, there exists a maximal incompressible subset $A\subset\tits X$ and $A$ is a compact convex spherical set \cite[Corollary~4.15]{GS_trans}.
Define $\mathcal C\subset\tits X$ to be the set of centers of maximal incompressible sets. 
By \cite[Proposition~5.2]{GS_trans}, there exists $\om\in \beta\Ga$ such that $T^\om$ maps the whole Tits boundary onto $A$.
In particular, $T^\om$ restricts to an isometry on every maximal incompressible set. Since $A$ is a convex spherical set, it follows that
$\mathcal C$ is a proper closed $\Ga$-invariant subset of $\geo X$. Hence, the action $\Ga\acts\geo X$ is not minimal.
\qed
\medskip

\proof[Proof of Theorem~\ref{thm_mainC}]
By Proposition~\ref{prop_reg_dich}, if there is no $\Ga$-periodic rank 1 geodesic, then the action $\Ga\acts\geo X$ is not minimal.
If this holds, then by \cite[Theorem~2.4]{ballmannbook}, the geodesic flow of $X$ cannot have a dense orbit.
Thus, Theorem~\ref{thm_mainB} completes the proof.
\qed

\proof[Proof of Theorem~\ref{thm_mainD}]
Let $X\cong X_0\times\ldots \times X_u$ be the deRham decomposition of $X$ \cite{FoLy_dR}.
By passing to a finite index subgroup, we may assume that $\Ga$ preserves the decomposition, namely every element $\ga\in\Ga$
has the form $\ga=(\ga_0,\ldots,\ga_u)$ with $\ga_i\in\Isom(X_i)$. For $i=1,\ldots, u$ define the groups $\Ga_i$ by projecting $\Ga$ to $X_i$. 
Then each $\Ga_i$ satisfies the duality condition in $X_i$. Moreover, 
since $\Ga$ acts cocompactly on $X$, every factor $X_i$  has a cocompact isometry group.
 By Proposition~\ref{prop_morse_rank}, the Tits boundary $\tits X\cong\tits X_0\circ\ldots\circ\tits X_u$ contains a regular point.
Thus, each factor $\tits X_i$ contains a regular point.
Hence the claim follows from Theorem~\ref{thm_mainC}.
\qed

\section{Tits alternative}\label{sec_tits}

Any subgroup of the isometry group of a  irreducible locally compact Euclidean building of dimension $n\geq 3$ is linear and therefore satisfies
Tits alternative. However, this argument fails for groups acting on 2-dimensional locally compact Euclidean building.
However, if a group $\Ga$ acts geometrically on a 2-dimensional locally compact Euclidean building, then it does satisfy Tits alternative 
\cite[Theorem~8.10]{BaBr_orbi}. The only ingredient in the proof is the fact that $\Ga$-periodic 2-flats are dense in the space of all 2-flats \cite[Theorem~8.9]{BaBr_orbi}.
We prepare for the proof of Theorem~\ref{thm_mainE} by establishing this fact for arbitrary buildings\footnote{This is certainly well-known but we were unable to
find a reference.}.

\bprop\label{prop_dense_tori}
Let $Y$ be a locally compact Euclidean building with a geometric group action $\Ga\acts Y$. Suppose that $\Ga$
satisfies the duality condition. Let $K$ be a compact subset of a flat $F\subset Y$. Then there exists a $\Ga$-periodic flat $F'$
which contains $K$.
\eprop

\proof
We may assume that $F$ is a maximal flat. Let $c\subset F$ be a regular geodesic with $o\in c$. Let $(\ga_k)$ be sequence in $\Ga$ such that 
$\ga_k^{\pm 1}o\to c(\pm\infty)$. For $k$ large enough, the piecewise geodesic $c_k$ determined by the points $\{\ga_k^jo\}_{j\in\Z}$
is a uniform Morse quasi-geodesic with arbitrary good Morse data \cite{KLP_morse}. By the Morse lemma for regular quasi-geodesics, $c_k$
lies in a uniform neighborhood of a unique maximal flat $F_k\subset Y$. Note that since the angles $\angle_o(\ga^{-1}_k o,\ga_k o)\to \pi$, and 
the action $\Ga\acts Y$ is properly discontinuous, the isometry $\ga_k$ is axial for large enough $k$ (compare \cite[Theorem~11]{Sw_cut}).
It follows that $F_k$ is a union of $\ga_k$-axes, since $\ga_k$ preserves $c_k$. Again, since $\Ga$ acts geometrically, it follows that $F_k$
is $\Ga$-periodic. For given $\eps>0$ we can find $k_\eps\in\N$ such that $c_k$ lies in the $r$-neighborhood of $F_k$ and  
$B_{\frac{1}{\eps}}(o)\cap c_k\subset N_\eps(c)$. Since $F$ and all $F_k$ are maximal flats, it follows that for large enough $k$ 
the compact set $K$ has to lie in $F_k$ as claimed. 
\qed

\proof[Proof of Theorem~\ref{thm_mainE}]
We may assume that $\Ga$ preserves the deRham decomposition of $X$.
By Theorem~\ref{thm_mainD} $X$ is a product of  a Euclidean space, rank 1 spaces, a Riemannian symmetric space and a Euclidean building.
If the deRham decomposition contains a rank 1 factor, then $\Ga$ contains a non-abelian free subgroup \cite[Theorem~3.5]{ballmannbook}.
The same is true if it contains a Riemannian symmetric space factor or a Euclidean building factor of dimension at least $3$.
This follows from Tits original free subgroup theorem since the corresponding groups are linear.
It remains to treat the case where $X$ is a product of a Euclidean space and irreducible Euclidean buildings of dimension 2:
\[X\cong \R^j\times Y_1\times\ldots\times Y_l\]
By Proposition~\ref{prop_dense_tori}, $\Ga$-periodic maximal flats are dense in the space of maximal flats.
Since $\Ga$ preserves the splitting, each $\Ga$-periodic maximal flat is a product of $\Ga_i$-periodic maximal flats in $Y_i$.
It follows that in each Euclidean building factor $Y_i$ the $\Ga_i$-periodic $2$-flats are dense in the space of $2$-flats.
This is the only ingredient in the proof of \cite[Theorem~8.10]{BaBr_orbi} which shows the existence of a non-abelian free subgroup
in case there is at least one Euclidean building factor.
\qed

\bibliographystyle{alpha}
\bibliography{rr_IV}

\Addresses

\end{document}